\documentclass[11pt,a4paper]{article}
\usepackage{epsf,epsfig,amsfonts,amsgen,amsmath,amstext,amsbsy,amsopn,amsthm,lineno}
\usepackage{color}

\newtheorem{theorem}{Theorem}

\newtheorem{false statement}{False statement}

\theoremstyle{definition}

\newtheorem{conjecture}{Conjecture}

\newtheorem{case}{Case}
\newtheorem{subcase}{Case}[case]

\begin{document}
\title{\bf\Large {A note on heavy cycles in weighted digraphs}\thanks{Supported  by NSFC (No.~10871158).}}

\date{}

\author{Binlong Li and Shenggui Zhang\thanks{Corresponding author.
E-mail address: sgzhang@nwpu.edu.cn (S. Zhang).}\\[2mm]
\small Department of Applied Mathematics,
\small Northwestern Polytechnical University,\\
\small Xi'an, Shaanxi 710072, P.R.~China\\}
\maketitle

\begin{abstract}
A weighted digraph is a digraph such that every arc is assigned a nonnegative number, called the weight of the arc. The weighted outdegree of a vertex $v$ in a weighted digraph $D$ is the sum of the weights of the arcs with $v$ as their tail, and the weight of a directed cycle $C$ in $D$ is the sum of the weights of the arcs of $C$. In this note we prove that if every vertex of a weighted digraph $D$ with order $n$ has weighted outdegree at least 1, then there exists a directed cycle in $D$ with weight at least $1/\log_2 n$. This proves a conjecture  of Bollob\'{a}s and Scott up to a constant factor.

\medskip

\smallskip
\end{abstract}

\section{Introduction}

We use Bondy and Murty \cite{Bondy_Murty} for terminology and
notation not defined here, and consider digraphs containing no
multiple arcs only.

Let $D$ be a digraph. The number of vertices and loops of $D$ are denoted by $n(D)$ and $r(D)$, respectively. We call $D$ a {\em
weighted digraph} if each arc $a$ of $D$ is assigned a nonnegative number
$w_D(a)$, called the {\em weight} of $a$. For a subdigraph $H$ of
$D$, $V(H)$ and $A(H)$ are used to denote the {set} of vertices and
arcs of $H$, respectively. The {\em weight} of $H$ is defined by
$$
w_D(H)=\sum_{a\in A(H)}w_D(a).
$$
For a vertex $v\in V(D)$, $N_H^+(v)$ denotes the set, and $d_H^+(v)$
the number, of vertices in $H$ to which there is an arc from $v$. We
define the {\em weighted outdegree} of $v$ in $H$ by
$$
d_H^{w+}(v)=\sum_{h\in N_H^+(v)}w_D(vh).
$$
When no confusion occurs, we will denote $w_D(a)$, $w_D(H)$,
$N_D^+(v)$, $d_D^+(v)$ and $d_D^{w+}(v)$ by $w(a)$, $w(H)$,
$N^+(v)$, $d^+(v)$ and $d^{w+}(v)$, respectively.

An unweighted digraph $D$ can be regarded as a weighted digraph in
which each arc $a$ is assigned weight $w(a)=1$. Thus, in an
unweighted digraph, $d^{w+}(v)=d^+(v)$ for every vertex $v$, and the
weight of a subdigraph is simply the number of its arcs.

A loopless digraph is one containing no loops. Let $D$ be a loopless
digraph such that every vertex of $D$ has outdegree at least $d$. {It is easy to see} that $D$ contains a {directed} cycle with length at least
$d+1$. For weighted digraphs, Bondy {\cite{Bondy}} conjectured that if every vertex
in a weighted loopless digraph has weighted outdegree at least 1,
then the digraph contains a {directed} cycle of weight at least 1. This conjecture was disproved by T. Spencer of Nebraska {(See \cite{Bollobas_Scott})}.

Bollob\'{a}s and Scott \cite{Bollobas_Scott} gave a lower bound
{on} the weight of heaviest directed cycles in a
weighted loopless digraph under the weighted outdegree condition.

\begin{theorem}[Bollob\'{a}s and Scott \cite{Bollobas_Scott}]
Let $D$ be a weighted loopless digraph with {$n\geq
2$} vertices. If $d^{w+}(v)\geq 1$ for every vertex $v\in V(D)$,
then $D$ contains a {directed} cycle $C$ such that
$w(C)\geq(24n)^{-1/3}$.
\end{theorem}

For an upper bound, Bollob\'{a}s and Scott
constructed a class of digraphs with minimum weighted outdegree at least
1 such that the maximum weight of cycles in these digraphs is at most
$c\log_2\log_2n/\log_2n$, where $c$ is a constant and $n$ is the
order of the digraph.

{As remarked} in \cite{Bollobas_Scott}, it seems
likely that $n^{-1/3}$ is much too small. Bollob\'{a}s and Scott
proposed the following conjecture.

\begin{conjecture}[Bollob\'{a}s and Scott \cite{Bollobas_Scott}]
Let $D$ be a weighted loopless digraph with {$n\geq
2$} vertices. If $d^{w+}(v)\geq1$ for every vertex $v\in V(D)$, then
$D$ contains a {directed} cycle $C$ such that $w(C)\geq 2/\log_2 n$.
\end{conjecture}

In this paper, we {prove the conjecture up to a
constant factor}.

\begin{theorem}
Let $D$ be a weighted loopless digraph with {$n\geq
2$} vertices. If $d^{w+}(v)\geq 1$ for every vertex $v\in V(D)$,
then $D$ contains a {directed} cycle $C$ such that $w(C)\geq
1/\log_2 n$.
\end{theorem}

{In fact, we can prove the following stronger assertion.}

\begin{theorem}
Let $D$ be a weighted digraph with {$n\geq 1$}
vertices and $r$ loops. If $d^{w+}(v)\geq 1$ for every vertex $v\in
V(D)$, then $D$ contains a {directed} cycle $C$ such that $w(C)\geq
1/\log_2 (n+r)$.
\end{theorem}

We postpone  the  proof of Theorem 3 to the next section.

\section{Proof of Theorem 3}

We use induction on $n$.

If $D$ {has} only one vertex, denote it by $v$. By $d^{w+}(v)\geq 1$,
we have $A(D)=\{vv\}$, $w(vv)\geq 1$ and $r=1$. Thus $C=vv$ is a
{directed} cycle with weight at least 1. The result is true. Now, we
suppose that $D$ has $n\geq 2$ vertices and $r$ loops.

\begin{case}
$D$ is not {strongly} connected.
\end{case}

Let $D'$ be a {strongly} connected component of $D$ such that there are
no arcs from $V(D')$ to $V(D)\backslash V(D')$. It is easy to know
that $d_{D'}^{w+}(v)=d_{D}^{w+}(v)\geq 1$ for all $v\in V(D')$. By
the induction hypothesis, there exists a {directed} cycle $C$ in $D'$
(and then, in $D$) such that
{$w(C)\geq {1}/{\log_2(n(D')+r(D'))}$}. Clearly $n(D')\leq n$ and
$r(D')\leq r$. Thus, we have {$w(C)\geq {1}/{\log_2(n+r)}$}, and complete the proof.

\begin{case}
$D$ is {strongly} connected.
\end{case}

\begin{subcase}
There exists a vertex $z$ such that $zz\notin A(D)$.
\end{subcase}

By the {strongly-connectedness} of $D$, there exists at least one arc
with head $z$. Let $y$ be a vertex such that $yz\in A(D)$ and
$w(yz)=\max\{w(vz): vz\in A(D)\}$. Consider the digraph $D'$ such
that {$V(D')=V(D)\backslash\{y\}$, $A(D')=A(D-y)\cup\{vz: vy\in
A(D)\}$,} and
$$
w_{D'}(uv)=\left\{
\begin{array}{ll}
  w_D(uy)+w_D(yz),  & \mbox{if\ } uy\in A(D) \mbox{\ and\ } v=z;\\
  w_D(uv),          & \mbox{otherwise}.
\end{array}
\right.
$$
Note that if $zy\in A(D)$, then $zz\in A(D')$, and
$w_{D'}(zz)=w_D(zy)+w_D(yz)$.

For every vertex $v\in V(D')$, its weighted outdegree is not less
than that in $D$. Thus, we have $d_{D'}^{w+}(v)\geq 1$ for all $v\in
V(D')$. By the induction hypothesis, there exists a {directed} cycle
$C'$ in $D'$ such that
{$w_{D'}(C')\geq {1}/{\log_2(n(D')+r(D'))}$}. {Since $n(D')=n-1$,
and $D'$ contains at most one loop more than} $D$, we have $r(D')\leq r+1$.
Thus {$w_{D'}(C')\geq {1}/{\log_2(n+r)}$}.

If $C'$ does not contain the vertex $z$, then it is also a {directed}
cycle in $D$ with the same weight. Otherwise, let $xz$ be the arc in
$C'$ with head $z$. If $xy\notin A(D)$, then $C'$ is also a {directed}
cycle in $D$ with the same weight. If $xy\in A(D)$, let
$C$ be the {directed} cycle obtained from $C'$ by {replacing the arc $xz$
with} the path $xyz$, then $C$ is a {directed} cycle in $D$ of weight
{$w_D(C)=w_{D'}(C')\geq {1}/{\log_2(n+r)}$.}

\begin{subcase}
For every $v\in V(D)$, $vv\in A(D)$.
\end{subcase}

In this case, $D$ has $r=n$ loops. And we need only prove that there
exists a {directed} cycle in $D$ with weight at least
{${1}/{\log_2(n+n)}={1}/{(1+\log_2n})$.}

If there exists a loop with weight at least {${1}/({1+\log_2 n})$},
then we complete the proof. So we assume that {every loop} of $D$
has weight less than {${1}/({1+\log_2 n})$}.

Let $D'$ be the digraph obtained from $D$ by deleting all the loops.
Then $D'$ has $n$ vertices and no loops, and for each vertex $v$ in
$V(D')$, {we have}
$$
d_{D'}^{w+}(v)\geq 1-\frac{1}{1+\log_2 n}=\frac{\log_2
n}{1+\log_2 n}.
$$

{It is easy to know that $D'$ is strongly
connected. Note that for every vertex $v\in V(D')$, $vv\notin
A(D')$. Using the conclusion of Case 2.1, we can obtained that}
there exists a {directed} cycle $C$ in $D'$ such that
$$
w_{D'}(C)\geq\frac{1}{\log_2
n}\frac{\log_2 n}{1+\log_2 n}=\frac{1}{1+\log_2 n},
$$
and $C$ is also a {directed} cycle in $D$ with the same weight.

The proof is complete.\hfill$\Box$

\end{document}